\documentclass[12pt,reqno]{amsart}

\usepackage{ifthen}

\newboolean{workmode}
\setboolean{workmode}{false}

\newboolean{sections}
\setboolean{sections}{true}

\usepackage[latin1]{inputenc} 
\usepackage[T1]{fontenc} 

\usepackage{amscd,graphicx}
\usepackage[all]{xy}
\usepackage{fullpage}
\usepackage{subsupscripts}

\usepackage[indent=0.5cm]{parskip}

	\AtBeginDocument{\addtocontents{toc}{\protect\setlength{\parskip}{0pt}}}
	\makeatletter
		\renewcommand*\l@subsection{\@tocline{2}{0pt}{2.75pc}{5pc}{}}
	\makeatother

\usepackage{amsmath}
\usepackage{amssymb}
\usepackage{amsthm}
\usepackage{thmtools} 
\usepackage{tikz-cd} 
\usepackage{rotating} 

\usepackage{xcolor}
\definecolor{jaw}{rgb}{0,.5,0}  
\definecolor{forestgreen}{rgb}{.2,.6,.2} 

\usepackage[colorlinks=true,linkcolor={red},citecolor={forestgreen},urlcolor={blue}]{hyperref}
\usepackage[capitalise,nameinlink,noabbrev]{cleveref}
	\crefname{subsection}{Subsection}{subsections}
	\crefdefaultlabelformat{#2\textbf{#1}#3} 

\newcommand{\ifwork}[1]{\ifthenelse{\boolean{workmode}}{#1}{}}
\newcommand{\comment}[1]{}
\newcommand{\mute}[1]{}
\newcommand{\printname}[1]{}

\ifwork{
	\renewcommand{\comment}[1]{{\marginpar{*}\ \scriptsize{#1}\ }}
	\renewcommand{\printname}[1]
	{{\color{brown}{\makebox[0pt]{\hspace{-1in}\raisebox{8pt}{\tiny #1}}}}}
}

\newcommand{\labell}[1]{\label{#1} \printname{#1}}

		\newcommand{\calD}{{\mathcal{D}}}

		\newcommand{\calX}{{\mathcal{X}}}
		\newcommand{\calY}{{\mathcal{Y}}}

		\newcommand{\QQ}{{\mathbb{Q}}}  
		\newcommand{\RR}{{\mathbb{R}}} 
		\newcommand{\R}{{\mathbb{R}}}  
		\renewcommand{\SS}{{\mathbb{S}}} 
		\newcommand{\ZZ}{{\mathbb{Z}}}  

		\newcommand{\coarse}{{\mathbf{Coarse}}} 
		\newcommand{\concr}{\mathbf{\hat{\kappa}}}
		\newcommand{\csheaves}{{\mathbf{CShf(\mfld)}}} 
		\newcommand{\diffeol}{\mathbf{Diffeol}} 
		\newcommand{\discr}{\mathbf{\hat{\pi}_0}} 
		\newcommand{\ext}{\operatorname{ext}} 
		\newcommand{\Hom}{{\operatorname{Hom}}} 
		\newcommand{\mfld}{{\mathbf{Man}}} 
		\newcommand{\open}{{\mathbf{Open}}} 
		\newcommand{\pizero}{\mathbf{{\pi}_0}} 
		\newcommand{\presheaves}{{\mathbf{PShf(\mfld)}}} 
		\newcommand{\res}{\operatorname{res}} 
		\newcommand{\set}{{\mathbf{Set}}} 
		\newcommand{\sheaves}{{\mathbf{Shf(\mfld)}}} 
		\newcommand{\stacks}{{\mathbf{St(\mfld)}}} 
		
		\newcommand{\G}{{\operatorname{G}}} 
		\newcommand{\SO}{{\operatorname{SO}}}  
		\newcommand{\U}{{\operatorname{U}}}  
		
		\newcommand{\CIN}{{C^\infty}}   
		\newcommand{\toto}{{~\rightrightarrows~}} 
		\newcommand{\rra}{{\rightrightarrows}}

		\usepackage{faktor}

		\makeatletter
		\DeclareRobustCommand*{\mfaktor}[3][]
		{
		   { \mathpalette{\mfaktor@impl@}{{#1}{#2}{#3}} }
		}
		\newcommand*{\mfaktor@impl@}[2]{\mfaktor@impl#1#2}
		\newcommand*{\mfaktor@impl}[4]{
		   \settoheight{\faktor@zaehlerhoehe}{\ensuremath{#1#2{#3}}}%
		   \settoheight{\faktor@nennerhoehe}{\ensuremath{#1#2{#4}}}%
		      \raisebox{-0.5\faktor@zaehlerhoehe}{\ensuremath{#1#2{#3}}}%
		      \mkern-4mu\diagdown\mkern-5mu%
		      \raisebox{0.5\faktor@nennerhoehe}{\ensuremath{#1#2{#4}}}%
		}
		\makeatother

	\newcommand{\ifsection}[2]{\ifthenelse{\boolean{sections}}{#1}{#2}} 

\theoremstyle{plain}

		\theoremstyle{plain}
		\ifsection{
			\newtheorem{theorem}{Theorem}[section] 
		}
		{
			\newtheorem{theorem}{Theorem} 
		}
		\newtheorem*{theorem*}{Theorem}
		\newtheorem{proposition}[theorem]{Proposition}
		\newtheorem{corollary}[theorem]{Corollary}
		\newtheorem*{conjecture*}{Conjecture}
		
		\newtheorem{lemma}[theorem]{Lemma}
		\newtheorem*{lemma*}{Lemma}

		\theoremstyle{definition}
		\newtheorem{definition}[theorem]{Definition}
		\newtheorem{example}[theorem]{Example}

		\newtheorem{notation}[theorem]{Notation}
		\newtheorem{remark}[theorem]{Remark}

		\def\eoe{\unskip\ \hglue0mm\hfill$/\!\!/$\smallskip\goodbreak} 

\author{Jordan Watts}
\address{Department of Mathematics, Central Michigan University, Mount Pleasant, Michigan, USA}
\email{jordan.watts@cmich.edu}

\author{Seth Wolbert*}
\address{Department of Mathematics and Statistics, Minnesota State University, Mankato, Minnesota, USA}
\email{seth.wolbert@mnsu.edu}

\title{Diffeological Coarse Moduli Spaces of Stacks over Manifolds}
\date{\today}

\begin{document}

\begin{abstract}In this paper, we consider diffeological spaces as stacks over the site of smooth manifolds, as well as the ``underlying'' diffeological space of any stack.  More precisely, we consider diffeological spaces as so-called concrete sheaves and show that the Grothendieck construction sending these sheaves to stacks has a left adjoint: the functor sending any stack to its diffeological coarse moduli space.  As an application, we restrict our attention to differentiable stacks and examine the geometry behind the coarse moduli space construction in terms of Lie groupoids and their principal bundles.  Within this context, we define a ``gerbe'', and show when a Lie groupoid is such a gerbe (or when a stack is represented by one).  Additionally, we define basic differential forms for stacks and confirm in the differentiable case that these agree (under certain conditions) with basic differential forms on a representative Lie groupoid.  These basic differentiable forms in turn match the diffeological forms on the orbit space.
\end{abstract}

\maketitle

\renewcommand*{\thefootnote}{\fnsymbol{footnote}}
\footnotetext{*This author acknowledges support from National Science Foundation grant DMS 08-38434 ``EMSW21-MCTP: Research Experience for Graduate Students''.}

\section{Introduction}\labell{s:intro}

The category of smooth manifolds lacks many structures that naturally arise in geometry.  Important examples include (critical) level sets of smooth real-valued functions, orbit spaces of Lie group actions, and function spaces, among many others.  Many categories extending the notion of smoothness have been introduced to account for this deficiency; for example, that of diffeological spaces, introduced by Souriau \cite{souriau} (see \cref{d:diffeology}).  Baez and Hoffnung in \cite{BH} show that the category of diffeological spaces forms a complete and co-complete quasi-topos; in particular, all subsets, quotients, and function spaces exist.

On the other hand, the language of stacks arose out of the need to better understand orbit and moduli spaces \cite{grothendieck}.  Indeed, the quotient stack associated to a compact Lie group action contains more information than the quotient topology, the ring of smooth functions on the orbit space in the sense of Schwarz \cite{CS,schwarz}, or even the associated diffeology.  For example, consider $\SO(n)$ acting on $\RR^n$ via rotations.  The orbit space is homeomorphic to $[0,\infty)$ for all $n$, and the rings of invariant smooth functions for each $n$ are all isomorphic.  The diffeology detects which $n$ we started with (this is the ``diffeological dimension'' at the singularity) \cite{iglesias}.  However, the orbit space obtained from the action of $\SO(2n)$ on $\RR^{2n}$ via rotations is diffeomorphic as a diffeological space to that obtained from the action of $\U(n)$ on the same space.  The corresponding stacks are not isomorphic in these examples, as the stabilizers at the origin are different, a property detected by the stacks.

The deficiencies of manifolds in modeling geometric structures has an effect on stacks over manifolds as well.  In particular, in the topological setting, a stack over the site of (compactly-generated) topological spaces admits a natural coarse moduli space \cite{noohi}.  This correspondence serves as a left adjoint to the ($2$-categorical) Yoneda embedding of topological spaces into stacks over topological spaces.  Unfortunately, in the case of stacks over manifolds, no such left adjoint exists: there are examples of stacks over manifolds for which no \emph{manifold} coarse moduli space exists \cite[Section 2]{hepworth:morse}. The purpose of this paper is to show that diffeology remedies this deficiency.

The first part this paper is mostly expository, casting results from category and topos theory known to experts into the languages of stacks over manifolds and diffeology, languages perhaps more familiar to differential geometers having a few references such as \cite{JY,johnstone}.  The second part explores the coarse moduli space of a stack over manifolds from a geometric perspective using principal bundles and diffeology, and is new.

As a bridge between the language of stacks and that of diffeology, we will utilize so-called concrete sheaves.  A concrete sheaf is a sheaf of sets that not only encodes an underlying set, but a family of maps that induces a diffeology on the set; see \cref{d:concrete sheaf}.  In fact, Baez and Hoffnung in \cite{BH} show that the category of diffeological spaces is equivalent to the category of concrete sheaves over open subsets of Euclidean spaces with smooth maps between them.  Moreover, Baez and Hoffnung show that the inclusion of concrete sheaves into all sheaves of sets has a left adjoint functor called concretization $\concr\colon\sheaves \to \csheaves$, which sends every sheaf of sets to its ``underlying'' concrete sheaf (or, equivalently, diffeology); see \cref{d:concretization}.

In this paper we go further, extending this adjunction to stacks over manifolds; note that we do not require the stacks to be differentiable.  Indeed, a straightforward computation shows that there is a functor $\coarse\colon\stacks \to \csheaves$ from the category of stacks over manifolds to the category of concrete sheaves, sending each stack to its diffeological coarse moduli space.  This is a left adjoint to $\int\colon\csheaves\to\stacks$, the Grothendieck construction taking concrete sheaves to stacks.  $\coarse$ is a 2-functor with respect to the 2-categorical structure of stacks and the degenerate 2-categorical structure of concrete sheaves (\emph{i.e.} that for which every 2-arrow the identity).  This 2-functor factors through a functor $\discr\colon\stacks\to\sheaves$ called discretization, which sends each stack to an ``underlying'' sheaf of sets, before being concretized in the sense described above.  Moreover, $\coarse$ has the useful property of being ``co-continuous'' in the sense that it sends (small) 2-colimits to (1-)colimits; see \cref{t:coarse co-continuous}.

Since manifolds form a full subcategory of diffeological spaces, they in turn form a full subcategory of concrete sheaves.  The Grothendieck construction embedding concrete sheaves into stacks over manifolds restricted to this full subcategory corresponds exactly to the ($2$-categorical) Yoneda embedding of manifolds into stacks.  In other words, for $\iota\colon\mfld\to\csheaves$ the aforementioned embedding of manifolds into concrete sheaves, $\int\colon\csheaves \to \stacks$ the Grothendieck construction, and $y\colon\mfld\to \stacks$ the Yoneda embedding, the diagram
\[\xymatrix{& \csheaves \ar[dr]^{\int} & \\ \mfld \ar[ur]^{\iota} \ar[rr]_y& & \stacks}\]
commutes.  For this reason, we may properly think of $\int\colon\csheaves \to \stacks$ as an extension of the Yoneda embedding.  As $\coarse$ and $\int$ form an adjunction, we recover a result similar to the topological coarse moduli space correspondence mentioned above (hence the name ``diffeological coarse moduli space'' for the object $\coarse(\mathcal{X})$ for any stack $\mathcal{X}$).  One may notice that $\iota$ here is the $1$-categorical Yoneda embedding of manifolds into sheaves over manifolds; thus $\int$ can be viewed as extending this embedding to its $2$-categorical counterpart.

Diffeology \emph{a priori} is a theory sitting in the realm of differential topology/geometry, not category theory.  
We therefore apply the categorical theory developed to differentiable stacks, obtaining geometric properties of their diffeological coarse moduli spaces.  If $G=(G_1\toto G_0)$ is a Lie groupoid representing a differentiable stack $\mathcal{X}$, we show that $\coarse(\mathcal{X})$ is (isomorphic to) the concrete sheaf associated to the quotient diffeology on $G_0/G_1$, the orbit space of the groupoid; see \cref{t:orbitspace}.  In the case where $\mathcal{X}=BG$,  the stack of right principal $G$-bundles, we show that discretization ($\discr$) identifies principal $G$-bundles over each manifold that are locally $G$-equivariantly bundle-diffeomorphic, and concretization ($\concr$) further identifies such bundles that are fiberwise $G$-equivariantly bundle-diffeomorphic; see \cref{p:loc eqvt diffeo,p:fiber eqvt diffeo}.

This paper is organized as follows.  \cref{s:diffeology} reviews diffeology and concrete sheaves, and extends their categorical equivalence to concrete sheaves over manifolds (\cref{l:concrete comparison}).  \cref{s:concrete} reviews that the inclusion of concrete sheaves into sheaves of sets over manifolds has a left adjoint functor (concretization), and that the inclusion of sheaves of sets into stacks over manifolds (also called the Grothendieck construction) has a left adjoint functor (discretization).  Together, these left adjoints compose to form $\coarse$ (see \cref{d:coarse} and \cref{t:coarse adj}).  We prove that it is co-continuous in \cref{t:coarse co-continuous}, in the sense that it takes 2-colimits (more prevalently used in stacks) to colimits.  As a first application, we restrict our attention to differentiable stacks in \cref{s:differentiable stacks}, obtaining the results mentioned above.  Here we define a ``gerbe'' (\cref{d:gerbe}) and show that the discretization and concretization of a differentiable stack are equal if and only if it is represented by such a gerbe (\cref{p:gerbiness}). In \cref{s:differential forms}, we consider basic differential forms as a second application, extending the results of \cite{KW, watts-groupoids} to certain differentiable stacks. 

For more details on diffeological spaces, see \cite{iglesias}.  For more details on stacks, see \cite{heinloth} and \cite{vistoli}, as well as \cite{lerman} for a description of the relationship between Lie groupoids and differentiable stacks. Coarse moduli spaces have been used in a variety of applications.  For example, for topological stacks, the \emph{topological coarse moduli space} is used in \cite{BGNX} to develop string topology for differentiable stacks over manifolds; in particular, the coarse moduli space is used to find ``hidden loops'' in the loop space of a stack.  In \cite{MT}, the authors use the topological coarse moduli space to extend the Kashiwara index formula from manifolds to orbifolds.  In \cite{hepworth:morse} and \cite{hepworth:vf}, the {\em underlying topological space of a stack} (similar to but fundamentally different from a coarse moduli space) is used to develop Morse theory and flows of vector fields for stacks over manifolds.  For more details on the topological coarse moduli space, see \cite{noohi}.  For this paper, we do not focus on topology; however, we do make note that the topology of the underlying topological space of a stack as defined by Hepworth in \cite{hepworth:morse} matches the so-called D-topology of the diffeological coarse moduli space (see \cref{r:topology}).

\subsection*{Prerequisites}\labell{ss:prerequisites}
This paper uses a lot of language from category theory.  Besides the basics, the reader is assumed to be familiar with the following.  We give an intuitive idea of what each of these objects is and pinpoint an actual reference for a definition or statement.  

\begin{itemize}

\item \labell{d:site} Given a category $C$, a \emph{coverage} on $C$ is a generalization of open covers from topology.  Essentially, it is a specification of the open covers of each object of $C$ satisfying some coherence conditions; see \cite[A2.1.9]{johnstone}.  A \emph{site} is a category equipped with a coverage.  (Note that we do not require the underlying category of a site to be small, whereas some references do.)  We will be dealing with the site of open subsets of Euclidean spaces with smooth maps as arrows (which is small), as well as the site of smooth manifolds (which is essentially small).  In both of these cases, the coverage is given by usual open covers by open sets.

\item \labell{d:sheaf} Given a category $C$, a \emph{presheaf of sets} on $C$ is a contravariant functor $C^{op}\to\set$.  If $C$ is a site, then a \emph{sheaf of sets} is a presheaf of sets on $C$ that obeys a ``sheaf condition'' with respect to the coverage on $C$; see \cite[A2.1.9]{johnstone}.  For example, given the site of manifolds, the sheaf $\Omega^k(\cdot)$ assigns to each manifold $M$ its set of $k$-forms $\Omega^k(M)$.  Presheaves of sets over $C$ form a category (and sheaves a subcategory) with arrows given by natural transformations.

\item \labell{d:sheafification} A presheaf of sets $P$ on a site $C$ has a unique \emph{sheafification}, a sheaf of sets $S$ over $C$ such that there is a natural transformation $\alpha\colon P\to S$ satisfying the following universal property: if $T$ is any sheaf of sets over $C$ and $\beta\colon P\to T$ a natural transformation, then $\beta$ factors through $\alpha$; see \cite[Definition 2.63 and Theorem 2.64]{vistoli}.  Note that sheafification is the left adjoint to the inclusion of sheaves of sets into presheaves of sets.

\item \labell{l:comparison} \textbf{Comparison Lemma:} Let $C$ be a locally small site, and $D$ a small dense subsite.  Then the restriction functor sending sheaves of sets over $C$ to sheaves of sets over $D$ is an equivalence of categories; see \cite[C2.2.3]{johnstone}.

\item \labell{d:cfg} Given a category $C$, a \emph{category fibered in groupoids (CFG)} over $C$ is a functor $\pi\colon D\to C$ whose fibers are groupoids, satisfying some conditions; see \cite[Definition 3.5]{vistoli}.  CFGs over $C$ form a $2$-category, where an arrow $F\colon(\pi\colon D\to C)\to(\pi'\colon D'\to C)$ is a functor from $D$ to $D'$ such that $\pi'\circ F=\pi$, and $2$-arrows are given by natural transformations. 

\item \labell{d:stack} A \emph{stack} over a site $C$ is a CFG that satisfies a ``descent condition'', which is a $2$-categorical version of the sheaf condition; see \cite[Definition 4.6]{vistoli}.  Stacks over $C$ form a sub-$2$-category of CFGs over $C$.

\item \labell{d:cat of el} Let $P\colon C^{op}\to\set$ be a presheaf of sets.  Then the \emph{category of elements} of $P$, denoted $\pi\colon \int P\to C$, is a CFG defined as follows.  The objects of $\int P$ are pairs $(x,c)$ where $c$ is an object in $C$ and $x\in P(c)$, and the arrows $\tilde{f}\colon(x,c)\to(x',c')$ are arrows $f\colon c\to c'$ in $C$ such that $P(f)(c')=c$.  Define $\pi$ to be the obvious projection.  Note that in this case, the resulting CFG is \emph{discrete}; that is, the arrows in each fiber are trivial, and so we think of these groupoids as just sets.  See \cite[page 203]{awodey} for more details.  It is not hard to show that $\int$ is a fully faithful functor from sheaves to CFGs, known as the \emph{Grothendieck construction}, sending a map of presheaves to the obvious map of CFGs.  Note that if $C$ is a site and $P$ a sheaf, then $\pi\colon \int P\to C$ is in fact a stack.

\end{itemize}

\subsection*{Acknowledgments}\labell{ss:acknowledgements}

The authors thank Eugene Lerman for many discussions regarding stacks, and introducing us to various categorical notions required to put the entire puzzle that is this paper together.  The authors also thank David M.\ Roberts for many useful suggestions, which resulted in a better paper.  The first author would also like to thank Patrick Iglesias-Zemmour and Yael Karshon for many illuminating diffeological examples, and Pierre Albin for discussions on the presentation.

\section{From Diffeology to Concrete Sheaves}\labell{s:diffeology}

In this section we review the definition of a diffeology and a concrete sheaf, and establish an equivalence of categories between diffeological spaces and concrete sheaves over manifolds (see \cref{c:diffeol equiv csheaves} and \cref{r:diffeol equiv csheaves}).  Our main resource for diffeology is the book by Iglesias-Zemmour \cite{iglesias}.  For more details on concrete sheaves, see \cite{BH}.

\begin{definition}
Let $\open$ be the category of open subsets of $\R^n$ for any $n$ and their smooth maps.  Let $\mfld$ be the category of manifolds.  Equipped with coverages generated by standard open covers, both $\open$ and $\mfld$ are sites.  Note that the singleton $\RR^0:=\{0\}$ is an object of $\open$ and $\mfld$.
\end{definition}

\begin{definition}[Diffeological Space]\labell{d:diffeology}
Let $X$ be a set.  A \emph{parametrization} $p\colon U\to X$ is a function from an object of $\open$ to $X$.  A \emph{diffeology} $\mathcal{D}$ on $X$ is a family of parametrizations satisfying the following three axioms:

\begin{itemize}
\item \textbf{(Covering Axiom)} All constant maps $p\colon U\to X$ are contained in $\mathcal{D}$.
\item \textbf{(Gluing Axiom)} If $p\colon U\to X$ is a parametrization such that there exist an open cover $\{U_\alpha\}_{\alpha\in A}$ of $U$ and a family $\{p_\alpha\colon U_\alpha\to X\}_{\alpha\in A}\subseteq\mathcal{D}$ of parametrizations that satisfy $p_\alpha=p|_{U_\alpha}$, then $p$ is contained in $\mathcal{D}$.
\item \textbf{(Smooth Compatibility Axiom)} If $p\colon U\to X$ is in $\mathcal{D}$ and $f\colon V\to U$ is a smooth map between open sets of Euclidean spaces, then $f^*p=p\circ f$ is also contained in $\mathcal{D}$.
\end{itemize}
The pair $(X,\mathcal{D})$ is called a \emph{diffeological space}, and the elements of $\mathcal{D}$ are called \emph{plots}.
\end{definition}

\begin{definition}[Smooth Map]\labell{d:smoothmap}
Let $(X,\mathcal{D}_X)$ and $(Y,\mathcal{D}_Y)$ be diffeological spaces, and let $F\colon X\to Y$ be a map.  Then $F$ is \emph{(diffeologically) smooth} if for every plot $p\in\mathcal{D}_X$, the composition $F\circ p$ is a plot in $\mathcal{D}_Y$.
\end{definition}

\begin{remark}\labell{r:diffeol is a cat}
Diffeological spaces with smooth maps between them form a category, denoted by $\diffeol$.  As mentioned in the introduction, this category is complete and cocomplete \cite{BH}, and contains the category of smooth manifolds with smooth maps $\mfld$ as a full subcategory.
\end{remark}

\begin{notation}
For $C$ a site, we will use the symbols $\mathbf{PShf}(C)$, $\mathbf{Shf}(C)$, $\mathbf{CFG}(C)$, and $\mathbf{St}(C)$ to denote the categories of presheaves of sets, sheaves of sets, categories fibered in groupoids, and stacks over $C$, respectively.
\end{notation}

\begin{definition}[Concrete Sheaves, Definition 4.10 of \cite{BH}]\labell{d:concrete sheaf}
Let $C$ be the site $\open$ or $\mfld$.  A presheaf $P\colon C^{op}\to\set$ is {\em concrete at $U\in C_0$} if the map ``\_'' sending each element $p\in P(U)$ to the function
\begin{equation}\labell{e:concrete}
\underline{p}\colon\Hom(*,U)\to P(*):u\mapsto u^*p
\end{equation}
is injective.  We say that $P$ is a {\em concrete presheaf} if it is concrete at each object $U$.  If in addition $P$ satisfies the sheaf axiom, then we call $P$ a \emph{concrete sheaf}.  We will refer to the category of concrete sheaves over $\open$ and $\mfld$ (with natural transformations as arrows) by $\mathbf{CShf}(\open)$ and $\csheaves$, respectively.
\end{definition}

\begin{lemma}[Baez-Hoffnung \cite{BH}, Prop. 4.15]\labell{l:baez-hoffnung}\labell{l:BH1}
There is an equivalence of categories from $\diffeol$ to $\mathbf{CShf}(\open)$.  This equivalence sends a diffeological space $X$ to the sheaf over $\open$ assigning to $U\in\open_0$ the set of all plots of $X$ with domain $U$.
\end{lemma}

The proof of the following lemma is a straightforward calculation and hence is omitted.

\begin{lemma}\labell{l:local lemma}
Let $C$ be the site $\open$ or $\mfld$ and $S\colon C^{op} \to \set$ a sheaf.  Let $U\in C_0$ and $\{U_i\to U\}$ a cover of $U$.  If $S$ is concrete at each $U_i$, then it is concrete at $U$.
\end{lemma}

\begin{lemma}[Concrete Comparison Lemma]\labell{l:concrete comparison}
The categories $\mathbf{CShf}(\open)$ and $\csheaves$ are equivalent.
\end{lemma}

\begin{proof}
Recall that the Comparison Lemma establishes an equivalence of categories
\begin{equation}\labell{comparisoneq}
  \xymatrix{\mathbf{Shf(\open)} \ar@<1ex>[r]^{\ext} & \sheaves \ar@<1ex>[l]^{\res}}
\end{equation}
where $\res(S)$ is the restriction of $S\colon\mfld^{op}\to \set$ to elements of $\open$ in $\mfld$ and $\ext$ takes $S'\colon\open^{op}\to \set$ to the Kan extension against the inclusion of $\open$ into $\mfld$.

The functor $\res$ takes $\csheaves$ to $\mathbf{CShf}(\open)$.  Since any $n$-dimensional manifold $M$ has an open cover by open subsets of $\RR^n$, it follows from \cref{l:local lemma} that the sheaf $S\colon\mfld^{op}\to \set$ must be concrete if $\res(S)$ is concrete.  So $\res$ takes {\em only} concrete sheaves of $\sheaves$ to concrete sheaves of $\mathbf{Shf(\open)}$.  Therefore, since $\res\circ\ext$ must take $\mathbf{CShf}(\open)$ to $\mathbf{CShf}(\open)$ (concreteness is preserved under isomorphism), it follows that $\ext$ takes $\mathbf{CShf}(\open)$ to $\csheaves$.  Thus the equivalence of categories (\cref{comparisoneq}) restricts to our desired equivalence between
$\mathbf{CShf}(\open)$ and $\csheaves$.
\end{proof}

\begin{corollary}\labell{c:diffeol equiv csheaves}
The categories $\diffeol$ and $\csheaves$ are equivalent.
\end{corollary}

\begin{remark}\labell{r:diffeol equiv csheaves}
The construction of a functor $I$ from $\diffeol$ to $\csheaves$ that is part of the equivalence in \cref{c:diffeol equiv csheaves} is given as follows.  Given a diffeological space $(X,\mathcal{D})$ and a manifold $M$ of dimension $m$, define a sheaf $I(X,\mathcal{D})$ over $\mfld$ by setting $I(X,\mathcal{D})(M)$ to be the set of all maps $p\colon M\to X$ such that there exist an open cover $\{U_\alpha\}$ of $M$ in which $U_\alpha\subseteq\RR^m$ for each $\alpha$ and plots $q_\alpha\colon U_\alpha\to X$ in $\mathcal{D}$ satisfying $p|_{U_\alpha}=q_\alpha$.  It follows that $I(X,\calD)(M)$ is the set of all diffeologically smooth maps $M\to X$; we will refer to these as \emph{plots} of $X$.  This is the concrete sheaf over $\mfld$ associated to $\mathcal{D}$.  A diffeologically smooth map $F$ is sent via $I$ to a map of sheaves sending a plot $p$ to $F\circ p$.
\end{remark}

\begin{remark}\labell{r:YonedaExt}
The fully faithful embedding of manifolds into diffeologies takes any manifold $M$ to the diffeological space built on the underlying set of $M$ with plots all smooth maps from open subsets of Euclidean space to $M$.  So the concrete sheaf $I(M)$ (for $I$ the functor of the previous remark) is just the sheaf of smooth functions to $M$: $I(M)(N)=C^\infty(N,M)$.  Write $\iota$ for the composition of the inclusion of the embedding of manifolds into diffeologies with $I$.  By the definition of the Yoneda embedding $y\colon \mfld \to \stacks$ (see, for instance, \cite[Section 8.2]{JY}), the following diagram commutes
\[\xymatrix{& \csheaves \ar[dr]^{\int} & \\ \mfld \ar[ur]^{\iota} \ar[rr]_y& & \stacks}\]
where $\int\colon\csheaves\to \stacks$ is the Grothendieck construction.  Thus, we may think of $\int$ as an extension of the Yoneda embedding.
\end{remark}

\section{From Concrete Sheaves to Sheaves to Stacks}\labell{s:concrete}

In this section we begin by reviewing concretization.  This is a functor introduced by Baez and Hoffnung in \cite{BH}, where they prove that it is left adjoint  to the inclusion $\csheaves\to\sheaves$.  We then define discretization following \cite{noohi}, and show that it is left adjoint to the Grothendieck construction.  These two functors compose to form the functor $\coarse$.  It follows that $\coarse$ is left adjoint to the restricted Grothendieck construction (\cref{t:coarse adj}) -- since $\csheaves$ is a $1$-category, we need only consider this left adjunction in a strict sense and not a $2$-categorical sense.  We similarly give the same treatment to $\discr$ in \cref{l:discr adj}.  Moreover, $\coarse$ is co-continuous, taking 2-colimits to colimits (\cref{t:coarse co-continuous}).

\begin{definition}[Concretization]\labell{d:concretization}
Let $P$ be a presheaf over $\mfld$.  Define $\kappa(P)\colon\mfld^{op}\to \set$ to be the presheaf assigning to each manifold $M$ the quotient set $P(M)/\!\sim$, where $p_1\sim p_2$ if $\underline{p_1}=\underline{p_2}$ (see \cref{e:concrete} for a definition of $\underline{p_i}$).  Define the \emph{concretization} of $P$, denoted by $\concr(P)\colon\mfld^{op} \to \set$, to be the sheafification of $\kappa(P)$.
\end{definition}

\begin{remark}\labell{r:concretization}
It follows from the definition of sheafification that $\concr(P)$ is a concrete sheaf.
\end{remark}

\begin{lemma}[Baez-Hoffnung, Lemma 5.20 of \cite{BH}]\labell{l:concr adj}
$\concr\colon\sheaves\to\csheaves$ is a left adjoint functor to the inclusion functor $\csheaves\hookrightarrow\sheaves$.
\end{lemma}

We'll see below that concretization can have interesting effects, especially on sheaves coming from differentiable stacks.  However, some sheaves of sets are practically destroyed by the process:
\begin{example}[$\concr$ Applied to $\Omega^k$]\labell{x:differential forms}
Let $\Omega^k\colon\mfld^{op}\to \set$ be the sheaf of differential $k$-forms.  If $k=0$, then $\concr(\Omega^k)$ is precisely equal to $\CIN(-,\RR)$, the Yoneda embedding of $\RR$ in $\sheaves$, or equivalently, $I(\RR)$ where $I$ is as in \cref{r:diffeol equiv csheaves}.  If $k>0$, then for every manifold $M$, we have $\concr(\Omega^k)(M)=\{0\}$. \eoe
\end{example}

\begin{definition}[Discretization]\labell{d:discretization}
Let $\mathcal{X}$ be a stack over $\mfld$.  Define $\pizero(\mathcal{X})$ to be the presheaf over $\mfld$ assigning to each manifold $M$ the set of isomorphism classes in $\mathcal{X}(M)$.  Define $\discr(\mathcal{X})$ to be the sheafification of $\pizero(\mathcal{X})$. For any manifold $M$ and $\xi\in\mathcal{X}(M)$, denote by $[\xi]$ the image of the isomorphism class of $\xi$ in $\discr(\mathcal{X})(M)$ via the natural map $\pizero(\mathcal{X})\to\discr(\mathcal{X})$.
\end{definition}

\begin{remark}\noindent
\begin{enumerate}
\item In general, sheafification is necessary in the definition of discretization; see \cref{x:BZ2}.
\item Note that for any stack $\mathcal{X}$, $\pizero(\mathcal{X})(*)=\discr(\mathcal{X})(*)$.
\end{enumerate}
\end{remark}

\begin{lemma}\labell{l:discr adj}
$\discr\colon\stacks\to\sheaves$ is left adjoint to the Grothendieck construction $\int\colon\sheaves\to\stacks$.
\end{lemma}

\begin{remark}
It is already known that $\int$ has a left adjoint functor, but this is typically described in much more generality (see \cite{awodey}).
\end{remark}

\begin{proof}
We will show $\int\colon\sheaves\to\stacks$ has the universal property of a right adjoint.  Fix a stack $\mathcal{X}\in\stacks$.  Then there exists a natural map of CFGs $\phi_\mathcal{X}\colon \mathcal{X} \to \int \pizero \mathcal{X}$ with $\phi_\mathcal{X}(\xi)=(M,[\xi])$ for each $\xi\in\calX(M).$

Suppose $S\colon\mfld^{op} \to \set$ is a sheaf of sets and $\psi\colon \mathcal{X} \to \int S$ is any map of stacks.  Since the fibers of $\int S$ are sets, $\psi$ must be constant on isomorphism classes of the fibers of $\mathcal{X}$.  It readily follows that $\psi$ factors uniquely through $\phi_\mathcal{X}$.

Since $\int$ is fully faithful, the unique map of CFGs $\bar{\psi}\colon\int \pizero \mathcal{X} \to \int S$ satisfying $\psi=\bar{\psi}\circ\phi_\mathcal{X}$ is equal to $\int$ applied to a unique map of presheaves $f\colon\pizero \mathcal{X} \to S$.  Furthermore, the map $f$ must factor uniquely through the sheafification $\eta_\mathcal{X}\colon\pizero \mathcal{X} \to \discr \mathcal{X}$.  In other words, there exists a unique map of sheaves $g\colon\discr \mathcal{X} \to S$ with $f=g\circ \eta$.  Thus, $\int(g)\colon\int\discr\mathcal{X}\to \int S$ is the unique map of stacks such that $\int(g)\circ (\int(\eta_{\mathcal{X}})\circ\phi_\mathcal{X}) = \psi$.

It follows $\int\colon\sheaves \to \stacks$ satisfies the universal property of a right adjoint with unit of adjunction $\mathcal{X}\mapsto\int(\eta_\mathcal{X})\circ \phi_\mathcal{X}$ and left adjoint $\discr$.
\end{proof}

\begin{remark}\labell{r:2-functor}
One can check that $\pizero\colon\stacks\to\presheaves$ sends $2$-commuting diagrams in $\stacks$ to strictly commuting diagrams in $\presheaves$.  It follows that $\discr$ is a $2$-functor, where $\presheaves$ inherits a $2$-category structure with trivial $2$-arrows.
\end{remark}

\begin{definition}[The Functor $\coarse$]\labell{d:coarse}
Define the \emph{coarse (2-)functor} $\coarse:\stacks\to\csheaves$ as the composition $\concr\circ\discr$.  For any stack $\mathcal{X}$, we will call the concrete sheaf $\coarse(\mathcal{X})$ {\em the diffeological coarse moduli space of $\mathcal{X}$} (see \cref{r:coarse moduli space}).
\end{definition}

\begin{remark}\labell{r:coarse}
A straightforward computation shows that one only needs to sheafify once when applying $\coarse$; indeed, we can define $\coarse$ to be the composition $\concr\circ\pizero$.
\end{remark}

\begin{theorem}[$\coarse$ is a Left Adjunction]\labell{t:coarse adj}
$\coarse$ is left adjoint to the Grothendieck construction restricted to $\csheaves$.
\end{theorem}

\begin{proof}
This follows directly from the composition of the left adjoint functors of \cref{l:concr adj} and \cref{l:discr adj}.
\end{proof}

\begin{remark}[Diffeological Coarse Moduli Space]\labell{r:coarse moduli space}
As noted in \cref{r:YonedaExt}, the functor $\int\colon\csheaves\to\stacks$ is an extension of the Yoneda embedding of manifolds into stacks and, as shown in \cref{t:coarse adj}, the functor $\coarse$ is a left adjoint to this extension.  This justifies our reference to $\coarse(\calX)$ as the diffeological coarse moduli space of the stack $\calX$, following the nomenclature of Hepworth (see \cite[Section 2]{hepworth:morse}).  We append the adjective ``diffeological'' to ``coarse moduli space'' as a standard coarse moduli space should be an object in the underlying site of a stack, but $\coarse(\mathcal{X})$ need not be a manifold in general. 
\end{remark}

\begin{theorem}[Coarse is Co-continuous]\labell{t:coarse co-continuous}
If a stack $\mathcal{X}$ is equivalent to a 2-colimit $\operatorname{colim}F$ in the 2-category of stacks, then $\coarse(\mathcal{X})$ is isomorphic to the 1-colimit $\operatorname{colim}(\coarse\circ F)$.
\end{theorem}

\begin{proof}
It follows from \cref{r:2-functor} that $\coarse$ takes $2$-commutative diagrams in $\stacks$ to strictly commuting diagrams in $\csheaves$.  Therefore, the $2$-cocones used to define colimits in $\stacks$ are sent to cocones in $\csheaves$. The proof now follows from a standard argument for co-continuity of left adjoints (for instance, see the dual argument in \cite[Theorem 1, page 118]{MacLane}).
\end{proof}

\begin{remark}\labell{r:generalizations}
Baez and Hoffnung generalize ``concreteness'' by defining what it means to be a concrete site (see Definition 20 of \cite{BH}).  In this way, concretization (\cref{d:concretization}) generalizes to any sheaf of sets over a site satisfying these concreteness conditions for which \cref{l:discr adj} is still valid.  Similarly, discretization (\cref{d:discretization}) generalizes to stacks over any site, as does \cref{l:discr adj}.  Therefore, \cref{t:coarse adj} generalizes as well, providing an adjunction to the Yoneda embedding of any concrete site $C$ into the category of stacks over $C$.  However, for this paper we are interested in diffeology, and so for the sake of brevity, we restrict our attention to the site $\mfld$. 
\end{remark}

\begin{remark}[Topological Considerations]\labell{r:topology}  Given a diffeological space $(X,\mathcal{D})$ the \emph{D-topology} of $X$ is defined to be the strongest topology on $X$ for which all plots in $\mathcal{D}$ are continuous.  That is, $A\subseteq X$ is open if and only if for every $p\in\mathcal{D}$, the set $p^{-1}(A)$ is open in the domain of $p$.

Let the diffeological coarse moduli space of the stack $\mathcal{X}$ be identified with the diffeological space $(X,\mathcal{D})$.  In particular, $\pizero(\mathcal{X}(*))$ is identified with the set $X$.  Then $A\subseteq X$ is open if and only if for every manifold $M$ and object $\xi\in\mathcal{X}(M)$, we have $\underline{[\xi]}^{-1}(A)\subseteq M$ is open.  Applying the Yoneda embedding $y$, we have that $A$ is open if and only if $\underline{[\varphi]}^{-1}(A)\subseteq M$ is open for every manifold $M$ and every map of stacks $\varphi\colon y(M)\to\mathcal{X}$.  This matches the topology of {\em the underlying space of a stack} used by Hepworth in \cite{hepworth:morse}.
\end{remark}

\section{Differentiable Stacks}\labell{s:differentiable stacks}

This section assumes that the reader is familiar with Lie groupoids, their actions and principal bundles, the localization of the $2$-category of Lie groupoids via bibundles and its equivalence to differentiable stacks over manifolds.  (One could also consider localizations \`a la Pronk \cite{pronk} using generalized morphisms, or \`a la Roberts \cite{roberts} using anafunctors; to keep in the spirit of principal bundles, however, we choose bibundles.) Our main resource for this theory is Lerman's survey on the topic, \cite{lerman}; we refer to this for definitions of bibundles, $BG$, etc. In this section, we show that the quotient diffeology on the orbit space $G_0/G_1$ of a Lie groupoid $G$ matches the diffeology associated to $\coarse(BG)$ (\cref{t:orbitspace}).  In the remainder of the section, we discuss the relationship between principal $G$-bundles and plots of the orbit space $G_0/G_1$, and an application to the ``relation space'' of a Lie groupoid (\cref{p:gerbiness}) that allows us to define a ``gerbe''.

The following is a special case of a folklore theorem well-known to experts; we refer to the proof given in \cite{CLW}.

\begin{lemma}\labell{l:BG 2-colimit}
Let $G=\left(G_1\underset{t}{\overset{s}{\toto}} G_0\right)$ be a Lie groupoid, $BG$ its associated stack, and $\calX$ a discrete stack (see the end of the Introduction for a definition; informally, this means that it is a sheaf of sets).  There is a bijection between maps of stacks $BG\to\calX$ and objects $\xi\in\calX(G_0)$ satisfying $s^*\xi=t^*\xi$.
\end{lemma}

\begin{proof}
By \cite[Proposition 6.6]{CLW}, there exists an equivalence of categories between cocycles $(\xi,f\colon s^*\xi\to t^*\xi)$ and maps of stacks $BG\to\calX$.  Since $\calX$ is discrete, $f$ is trivial and $\Hom(BG,\calX)$ has trivial arrows.  Thus, the equivalence reduces to a bijection.
\end{proof}

\begin{theorem}[The Orbit Space of a Lie Groupoid]\labell{t:orbitspace}
Let $\mathcal{X}$ be a differentiable stack.  Let $G=(G_1\toto G_0)$ be a Lie groupoid representing $\mathcal{X}$ (that is, $\mathcal{X}\cong BG$).  Then $\coarse(\mathcal{X})$ is isomorphic to the concrete sheaf associated to the quotient diffeology on the orbit space $G_0/G_1$ of $G$.
\end{theorem}

\begin{remark}
In particular, the above theorem shows that the diffeology associated to the diffeological coarse moduli space $\coarse(BG)$ is not only independent of the choice of an atlas for $BG$, but only dependent on the equivalence class of $BG$ in the category of \emph{all} stacks over manifolds.
\end{remark}

\begin{proof}
The quotient map $\pi\colon G_0\to G_0/G_1$ corresponds to the associated plot in $\int I(G_0/G_1)(G_0)$ via the $2$-Yoneda Lemma (see \cite[Section 8.2]{JY}), where $I$ is the functor of \cref{r:diffeol equiv csheaves}.  By \cref{l:BG 2-colimit}, there is a corresponding map $\theta\in\Hom(BG,\int I(G_0/G_1))$.  Applying $\coarse$ and \cref{c:diffeol equiv csheaves}, and using the universal property of diffeological quotients, we obtain an isomorphism $\coarse(BG)\cong I(G_0/G_1)$.
\end{proof}

Compare the following corollary with \cite{karshon-miyamoto} and \cite{watts-groupoids}.

\begin{corollary}\labell{c:orbitspace}
Given differentiable stacks $\calX$ and $\calY$ with a map of stacks $F\colon\calX\to\calY$, the map $\coarse(F)$ corresponds to a smooth map between the corresponding quotient diffeologies of representative Lie groupoid orbit spaces.  Moreover, if $F$ is an equivalence (and hence Lie groupoids representing $\calX$ and $\calY$ are Morita equivalent), then $\coarse(F)$ is a diffeomorphism.
\end{corollary}

Let $G$ be a Lie groupoid, and let $BG$ be its corresponding differentiable stack.  Given principal $G$-bundles $P$ and $Q$ over a manifold $B$, we have that $P$ and $Q$ are isomorphic (\emph{i.e.}\ $G$-equivariantly bundle-diffeomorphic) if and only if they represent the same element in $\pizero(BG)(B)$.  This equivalence no longer holds upon sheafification of $\pizero(BG)$ to $\discr(BG)$ (see \cref{x:BZ2}).

\begin{definition}[Locally Isomorphic Bundles]\labell{d:loc eqvt diffeo}
Let $G$ be a Lie groupoid, and let $P$ and $Q$ be principal $G$-bundles over a manifold $B$.  Then $P$ and $Q$ are \emph{locally isomorphic} if there exists an open cover $\{i_\alpha\colon B_\alpha\to B\}$ of $B$ such that for each $\alpha$, the principal $G$-bundles $i_\alpha^*P$ and $i_\alpha^*Q$ are isomorphic.
\end{definition}

\begin{proposition}[Discretization for Differentiable Stacks]\labell{p:loc eqvt diffeo}
Let $G$ be a Lie groupoid, and let $P$ and $Q$ be principal $G$-bundles over a manifold $B$.  Then, $P$ and $Q$ are locally isomorphic if and only if $[P]=[Q]$ in $\discr(BG)(B)$.
\end{proposition}

\begin{proof}
$P$ and $Q$ are locally isomorphic if and only if there exists an open cover $\{i_\alpha\colon B_\alpha\to B\}$ of $B$ such that for each $\alpha$, there is an isomorphism $i_\alpha^*P \cong i_\alpha^*Q$.  This is equivalent to the equality of the isomorphism class of $i_\alpha^*P$ with the isomorphism class of $i_\alpha^*Q$ in $\pizero(BG)(B_\alpha)$, for each $\alpha$.  By definition of sheafification, this is equivalent to $[P]=[Q]$.
\end{proof}

The following example illustrates the difference between $\pizero(BG)$ and $\discr(BG)$.

\begin{example}[$B\ZZ_2$]\labell{x:BZ2}
Let $\mathcal{X}$ be the differentiable stack $B\ZZ_2$, where $\ZZ_2=\ZZ/2\ZZ$.  Then $B\ZZ_2(\SS^1)$ contains two isomorphism classes of bundles: one represented by the trivial $\ZZ_2$-bundle which we denote by $P$, and the other represented by the boundary of the M\"obius band which we denote by $Q$.  While $P$ and $Q$ do not represent the same elements in $\pizero(BG)$, by \cref{p:loc eqvt diffeo} we know that $[P]=[Q]$ in $\discr(BG)$.
\eoe
\end{example}

\begin{remark}
The above example is not unique; given {\em any} Lie group $\Gamma$, local triviality of principal $\Gamma$ bundles ensures that all such bundles are locally isomorphic and hence $\discr(B\Gamma)(M)$ is a singleton set for each $M\in\mfld_0$.
\end{remark}

Given a Lie groupoid $G$ and two principal $G$-bundles $P$ and $Q$, we give conditions under which $[P]$ and $[Q]$ yield the same objects after applying the functor $\kappa$ to $\discr (BG)$.

\begin{definition}[Fiberwise Isomorphic Bundles]\labell{d:fiber eqvt diffeo}
Let $G$ be a Lie groupoid, and let $\rho_P\colon P\to B$ and $\rho_Q\colon Q\to B$ be principal $G$-bundles over a manifold $B$.  Then $P$ and $Q$ are \emph{fiberwise isomorphic} if for every $b\in B$, there exists a $G$-equivariant diffeomorphism from $\rho_P^{-1}(b)$ to $\rho_Q^{-1}(b)$.
\end{definition}

\begin{proposition}[Concretization for Differentiable Stacks]\labell{p:fiber eqvt diffeo}
Let $G$ be a Lie groupoid, and let $\rho_P\colon P\to B$ and $\rho_Q\colon Q\to B$ be principal $G$-bundles over a manifold $B$.  Then, $P$ and $Q$ are fiberwise isomorphic if and only if $P$ and $Q$ descend to the same elements in $\coarse(BG)(B)$.
\end{proposition}

\begin{proof}
$P$ and $Q$ descend to the same elements of $\kappa(\discr(BG))(B)$ if and only if $\underline{[P]}$ and $\underline{[Q]}$ are the same.  But this is equivalent to $[b^*P]=[b^*Q]$ for each $b\in B$. That is, for each $b\in B$, there is a $G$-equivariant diffeomorphism between $b^*P$ and $b^*Q$.  But these two objects are exactly the fibers $\rho_P^{-1}(b)$ and $\rho_Q^{-1}(b)$, respectively.  The result now follows from the fact that any concrete presheaf $P$ is \emph{separated} \cite[Lemma 5.21]{BH}; \emph{i.e.}\ its sheafification induces an injection $\kappa(P)(M)\to\concr(P)(M)$ for each manifold $M$.
\end{proof}

Not every pair of fiberwise isomorphic principal bundles are locally isomorphic, as illustrated by the following example.

\begin{example}[Reflection in $\RR$]\labell{x:fiber eqvt diffeo}
Let $G$ be the action groupoid $(\ZZ_2\times\RR)\toto\RR$ obtained from the action of $\ZZ_2$ on $\RR$ sending $(\pm1,x)$ to $\pm x$.  Consider the following two plots $p_1,p_2\colon\RR\to\RR$ defined as:

$$ p_1(\tau)=\begin{cases}
-e^{-1/\tau^2} & \text{ if $\tau<0$,}\\
0 & \text{ if $\tau=0$,}\\
e^{-1/\tau^2} & \text{ if $\tau>0$,}\\
\end{cases}$$
$$ p_2(\tau)=\begin{cases}
-e^{-1/\tau^2} & \text{if $\tau\neq 0$,}\\
0 & \text{ if $\tau=0$.}\\
\end{cases}$$

Computing the pullback bundles $P_i:=p_i^*t$ ($i=1,2$), where $t\colon G_1\to G_0$ is the unit bundle (in fact, this is a trivial $\ZZ_2$-bundle), we have
\begin{align*}
P_1=&~\{(\tau,a,x)\in\RR\times\ZZ_2\times\RR~|~p_1(\tau)=a\cdot x\}\\
=&~\{(\tau,1,-e^{-1/\tau^2})~|~\tau<0\}\cup\{(0,\pm1,0)\}\cup\{(\tau,1,e^{-1/\tau^2})~|~\tau>0\}\\
&~\cup\{(\tau,-1,e^{-1/\tau^2})~|~\tau<0\}\cup\{(\tau,-1,-e^{-1/\tau^2})~|~\tau>0\},
\end{align*}
and similarly,
\begin{align*}
P_2=&~\{(\tau,1,-e^{-1/\tau^2})~|~\tau<0\}\cup\{(0,\pm1,0)\}\cup\{(\tau,1,-e^{-1/\tau^2})~|~\tau>0\}\\
&~\cup\{(\tau,-1,e^{-1/\tau^2})~|~\tau<0\}\cup\{(\tau,-1,e^{-1/\tau^2})~|~\tau>0\}.
\end{align*}

If $\phi\colon P_1\to P_2$ is a bundle-diffeomorphism, then in particular it must send $\tau$ to $\tau$; that is, $\phi$ must fix the first argument of each triple in $P_1$.  The anchor map $\alpha_i\colon P_i\to \RR$ is the projection to the third argument of each triple, and so if $\phi$ is equivariant, then it must fix the third argument of each triple in $P_1$.  Thus, $\phi$ must be the identity map when restricted to fibers for which $\tau<0$, and it switches connected components of the fibers for $\tau>0$.  Hence, $\phi$ cannot be continuous at $\tau=0$, and so we conclude that there is no equivariant diffeomorphism of bundles; in fact, not even a local one exists.  However, it is clear that the fibers of $P_1$ and $P_2$ are equivariantly diffeomorphic, and so by \cref{p:fiber eqvt diffeo} we have $\underline{[P_1]}=\underline{[P_2]}$.  This corresponds to the fact that $\pi_G\circ p_1=\pi_G\circ p_2$ where $\pi_G\colon G_0\to G_0/G_1$ is the quotient map (see \cref{p:plots}). \eoe
\end{example}

In light of \cref{r:diffeol equiv csheaves}, we show in the following proposition how elements of $\coarse(BG)$ and plots of the corresponding diffeology are related.

\begin{proposition}[From Principal Bundles to Plots]\labell{p:plots}
Let $G$ be a Lie groupoid, and let $\rho\colon P\to B$ be a principal $G$-bundle, where $B$ is an open subset of some Euclidean space.  Then, $P$ induces a plot $p\colon B\to G_0/G_1$ defined as 
\begin{equation}\labell{e:plots}
p(b)=\pi_G\circ\alpha\circ\sigma(b)
\end{equation}
where $\pi_G\colon G_0\to G_0/G_1$ is the quotient map, $\alpha\colon P\to G_0$ is the anchor map, and $\sigma$ is a local section of $\rho$ about $b$.  Moreover, every plot in the quotient diffeology locally arises in this way.  In particular, $\pizero(BG(*))$ is in bijection with the underlying set of $G_0/G_1$.   Finally, two principal $G$-bundles are fiberwise isomorphic if and only if the corresponding plots are equal.
\end{proposition}

\begin{remark}\labell{r:plots}
The statement that $\pizero(BG(*))$ is in bijection with $G_0/G_1$ is not new, but we emphasize it due to it being crucial to understanding the diffeological coarse moduli space of a differentiable stack from the standpoint of diffeology.  Also, if $B$ is a manifold, but not an open subset of some Euclidean space, then the proposition still holds where $p$ is a plot in the sense that it is contained in $I(\mathcal{D})(B)$ where $\mathcal{D}$ is the quotient diffeology on $G_0/G_1$ and $I$ is as given in \cref{r:diffeol equiv csheaves}.
\end{remark}

\begin{proof}
Suppose that $\sigma$ and $\sigma'$ are two local sections of $\rho$ about $b$.  There exists $g\in G_1$ such that $\sigma(b)\cdot g=\sigma'(b)$.  Consequently, $\alpha(\sigma(b))$ and $\alpha(\sigma'(b))$ lie in the same $G$-orbit, and so descend to the same point via $\pi_G$.  The map $p$ is thus independent of the local section chosen, and so is well-defined.  Moreover, $p$ is diffeologically smooth, and so is a plot.  To show that $p$ and $\underline{[P]}$ are the same, consider the bundle $b^*P$, which is the fiber of $P$ over $b$.  The anchor map restricted to $b^*P$ has an image contained in exactly one $G$-orbit of $G_0$, which by definition of $G_0/G_1$ is exactly $p(b)$.  But $[b^*P]=\underline{[P]}(b)$, and this proves the claim.

On the other hand, given any plot $p\colon U\to G_0/G_1$ in the quotient diffeology, there exist an open cover $\{U_\alpha\}$ of $U$ and a family $\{q_\alpha\colon U_\alpha\to G_0\}$ of smooth maps such that $p|_{U_\alpha}=\pi_G\circ q_\alpha$.  Taking the pullback bundles $q_\alpha^*(t\colon G_1\to G_0)$ proves the claim.

The last statement follows immediately from \cref{e:plots}. 
\end{proof}

Given the above discussion, a natural question arises: for $BG$ a differentiable stack, when is $\discr(BG)$ a concrete sheaf; \emph{i.e.}\  when does $\coarse(BG)=\discr(BG)$?  It turns out that this condition is equivalent to asking that $BG$ be a differentiable gerbe over its quotient space.  To define this, we first need to define subductions, the natural generalization of a surjective submersion to the category of diffeological spaces.

\begin{definition}\labell{d:subduction}
A map $F\colon X\to Y$ of diffeological spaces is a \emph{subduction} if it is a smooth surjection such that for every plot $p\colon U\to Y$ and every $u\in U$ there is an open neighborhood $V$ of $u$ and a plot $q\colon V\to X$ such that $p|_V=F\circ q$.
\end{definition}

\begin{example}\labell{x:quotient}
Given a diffeological space $X$ and an equivalent relation $\sim$, the quotient map $X\to X/\!\sim$ is a subduction.  In particular, given a Lie groupoid $G_1\toto G_0$, the quotient map to the orbit space $G_0\to G_0/G_1$ is a subduction.
\eoe
\end{example}

\begin{remark}\labell{r:relations}
Suppose that $G=(G_1\to G_0)$ is a Lie groupoid and $X$ is a diffeological space.  Then by identifying $X$ with the trivial groupoid $X\toto X$, we may treat $X$ as a \emph{diffeological groupoid}; that is, a groupoid in which the object and arrow spaces are diffeological spaces, and all structure maps are smooth. 
 With this in mind, suppose that $\pi:G\to X$ is smooth functor (i.e., a functor which is smooth on objects and morphisms).  As a slight abuse of notation, we will also use $\pi$ to refer to the associated map $\pi\colon G_0\to X$.  Via the basic properties of functors and pullbacks, such a functor $\pi\colon G\to X$ induces a unique smooth map
\[(s,t)\colon G_1\to G_0\times_XG_0,\quad (s,t)(g):=(s(g),t(g))\]
\end{remark}

As a particular example of \cref{r:relations}, note that the quotient map $G_0\to G_0/G_1$ discussed in \cref{x:quotient} induces a smooth functor $G\to G_0/G_1$ which, due to its importance for what follows, we define separately with its own special notation.

\begin{definition}
Let $G=(G_1\toto G_0)$ be a Lie groupoid.  Then the {\em relation space of $G$} is the diffeological pullback
\[\xymatrix{R \ar[r] \ar[d] & G_0 \ar[d]^{\pi} \\G_0\ar[r]_-{\pi} & G_0/G_1}\]
for $\pi\colon G_0\to G_0/G_1$ the diffeological quotient map.
\end{definition}

We can now define what it means to be a gerbe over a diffeological space.  Note this is an extension of \cite[Definition 6.1]{RV}, where we have replaced surjective submerions with subductions.    

\begin{definition}\labell{d:gerbeoverX}
Let $G=(G_1\to G_0)$ be a Lie groupoid and $X$ be a diffeological space. A smooth functor $\pi\colon G\to X$ is a \emph{gerbe over $X$} if the maps $\pi\colon G_0\to X$ and $(s,t)\colon G_1\to G_0\times_XG_0$ are both subductions.
\end{definition}

While \cref{d:gerbeoverX} is a ``correct'' definition from the perspective of thinking about gerbes as objects over a base space, it turns out that a smooth functor $\pi\colon G\to X$ can only be (up to diffeomorphism) a gerbe over the quotient space $G_0/G_1$ of $G$.  Indeed, we have the following:

\begin{proposition}
Let $G=(G_1\to G_0)$ be a Lie groupoid, $\pi\colon G\to X$ be a smooth functor, and $q\colon G_0\to G_0/G_1$ be the diffeological quotient map.  If $\pi\colon G_0\to X$ is a subduction and $(s,t)\colon G_1\to G_0\times_X G_0$ is surjective, then $\pi$ is also a diffeological quotient map with respect to the relation on $G_0$ induced by $G_1$; that is, there exists a unique diffeomorphism $f\colon G_0/G_1\to X$ satisfying $f\circ q=\pi$.

In particular, if $\pi\colon G\to X$ is a gerbe over $X$, then $\pi\colon G_0\to X$ is a diffeological quotient map with respect to the relation on $G_0$ induced by $G_1$.
\end{proposition}

\begin{proof}
First, note that $\pi$ being a functor means that $\pi$ is invariant under the relation on $G_0$ induced by $G_1$.  Thus, via the universal property of quotient spaces, there exists a unique smooth map $f\colon G_0/G_1\to X$ for which $f\circ q = \pi$.

Next, note that $f$ is bijective: surjectivity follows from the surjectivity of $\pi$.  For injectivity, suppose that $q(x),q(y)\in G_0/G_1$ satisfy $f(q(x))=f(q(y))$.  Then we have $\pi(x)=\pi(y)$ and thus $(x,y)\in G_0\times_X G_0$.  Therefore, by the surjectivity of $(s,t)\colon G_1\to G_0\times_X G_0$, there exists $g\in G_1$ with $s(g)=x$ and $t(g)=y$ and thus $q(x)=q(y)$.

Finally, note that $f$ is a subduction: as $\pi$ is a subduction, any plot $p\colon U\to X$ locally lifts to a plot of $G_0$ which, after composing with $q$, descends to a local lift of $p$ to a plot of $G_0/G_1$.  Therefore, $f$ is a diffeomorphism (as is any injective subduction, see Article 1.49 of \cite{iglesias}).
\end{proof}

Thus, by the above proposition, a gerbe $\pi\colon G\to X\cong G_0/G_1$ over a diffeological space $X$ can be thought of as a property intrinsic to the groupoid itself and hence one may simply call the groupoid $G$ a gerbe and suppress the reference to the base space $X$.

\begin{definition}\labell{d:gerbe}
A Lie groupoid $G=(G_1\toto G_0)$ is a \emph{gerbe} if the map $(s,t)\colon G_1\to R$ is a subduction, where $R$ is the relation space of $G$.  Call a stack $\mathcal{X}$ a \emph{differentiable gerbe} if it is presented by a gerbe (i.e., if it is isomorphic to $BG$ as a stack). 
\end{definition}

\begin{example}[Irrational Torus]\labell{x:gerbe}
Fix $\alpha\in\RR\smallsetminus\QQ$.  Consider the action groupoid of the Lie group action of $\ZZ^2$ on $\RR$ given by $(m,n)\cdot x:= x+m+n\alpha$.  The relation space $R$ is a dense subset of $\RR^2$ (a dense collection of translates of the line $y=x$), and $(s,t)\colon \ZZ^2\times\RR\to R$ is a subduction.  Note that $T_\alpha:=\RR/\ZZ^2$ is \emph{not} a manifold, and so our concept of gerbe indeed generalizes \cite[Definition 6.1]{RV}.
\eoe
\end{example}

Finally, we will show that we can tie the property of being a differentiable gerbe to our functor $\coarse$.  To do so, we will need the following lemma on the local nature of principal $G$-bundles and their morphisms.

\begin{lemma}\labell{l:local isomorphisms are n t}
Let $G=(G_1\rra G_0)$ be a Lie groupoid and let $f_i\colon M\to G_0$ ($i=1,2$) be two smooth functions.  Then the pullback bundles $P_i:=f_i^*G_1$ are isomorphic if and only if there exists a smooth function $\varphi\colon M\to G_1$ satisfying $s\circ \varphi =f_1$ and $t\circ \varphi=f_2$.
\end{lemma}

\begin{proof}
Given a function $\varphi\colon M\to G_1$ satisfying $s\circ\varphi = f_1$ and $t\circ\varphi = f_2$, we have the isomorphism
\[\tilde{\varphi}\colon P_1\to P_2,\;(p,g)\mapsto (p,\varphi(p)\cdot g), \text{ for }(p,g)\in P_1=M\times_{f_1,G_0,t}G_1.\]

Conversely, recall that for any principal $G$-bundle $P$ over $M$, there exists a smooth division map 
\[d\colon P\times_{M} P \to G_1,\; (p,q)\mapsto g,\text{ for }g\text{ satisfying }p\cdot g = q\]
(this is a consequence of the requirement that the map $P\times_{G_0}G_1 \to P\times_{M}P$ sending $(p,g)$ to $(p,p\cdot g)$ be a diffeomorphism).  

Let $\psi\colon P_1\to P_2$ be an isomorphism of principal $G$-bundles, let $d$ designate the division map for $P_2$, and denote by $u\colon G_0\to G_1$ the unit map.  Let $\tilde{f}_i\colon M\to P_i$ be the maps $x\mapsto (x,u(f_i(x)))$.  Then 
\[\varphi\colon M\to G_1\colon x\mapsto d\left(\tilde{f}_2(x),\psi(\tilde{f}_1(x))\right) \]
yields the required map with $s\circ \varphi=f_1$ and $t\circ \varphi=f_2$.
\end{proof}

\begin{proposition}\labell{p:gerbiness}
Let $G$ be a Lie groupoid.  Then $\discr(BG)$ is a concrete sheaf of sets if and only if $G$ is a gerbe.
\end{proposition}

\begin{proof}
By \cref{p:loc eqvt diffeo}, $\discr(BG)$ is the sheaf of local isomorphism classes of principal $G$-bundles.  By \cref{p:fiber eqvt diffeo}, we have that any two principal $G$-bundles descend to the same plot of $\coarse(BG)=\concr(\discr(BG))$ exactly when they are fiberwise isomorphic.  Since $\discr(BG)$ is concrete exactly when $\concr(\discr(BG))= \discr(BG)$, it is enough to prove that the map $(s,t)$ is a subduction exactly when any two principal $G$-bundles which are fiberwise isomorphic are locally isomorphic.  

By definition, $(s,t)$ is a subduction if any plot $p$ of $R$ locally lifts to plots of $G_1$.  By definition of $R$, this is equivalent to the following: given a pair of plots $p_i\colon M\to G_0$ ($i=1,2$) satisfying $\pi\circ p_1=\pi\circ p_2$, the plot $(p_1,p_2)$ locally lifts to plots of $G_1$.  Let $\pi\colon G_0\to G_0/G_1$ be the quotient map.

Suppose $(s,t)$ is a subduction.  Let $P_i$ ($i=1,2$) be fiberwise isomorphic principal $G$-bundles.  There exist an open cover $\{U_\alpha\}$ of $M$ and for each $\alpha$ and $i=1,2$, a plot $p_i^\alpha\colon U_\alpha\to G_0$ such that $P_i|_{U_\alpha}\cong(p_i^\alpha)^*G_1$.  Fix $\alpha$.  Then $(p_1^\alpha)^*G_1$ is fiberwise isomorphic to $(p_2^\alpha)^*G_1$.  By \cref{p:plots}, $\pi\circ p_1^\alpha=\pi\circ p_2^\alpha$.  Thus, $(p_1^\alpha,p_2^\alpha)$ locally lifts to $G_1$.  By \cref{l:local isomorphisms are n t}, these local lifts to $G_1$ yield local isomorphisms from $P_1|_{U_\alpha}$ to $P_2|_{U_\alpha}$.  Therefore $P_1$ and $P_2$ are locally isomorphic.

Conversely, suppose any fiberwise isomorphic principal $G$-bundles are locally isomorphic.  Let $p_i\colon M\to G_0$ ($i=1,2$) be plots satisfying $\pi\circ p_1=\pi\circ p_2$.  By \cref{p:plots}, $p_1^*G_1$ is fiberwise isomorphic to $p_2^*G_2$.  By hypothesis, these are locally isomorphic; that is, by \cref{l:local isomorphisms are n t}, there is an open cover $\{U_\alpha\}$ of $M$ and plots $g_\alpha\colon U_\alpha\to G_1$ such that $s\circ g_\alpha=p_1|_{U_\alpha}$ and $t\circ g_\alpha=p_2|_{U_\alpha}$, for each $\alpha$.  So $(s,t)\circ g_\alpha=(p_1,p_2)|_{U_\alpha}$ for each $\alpha$, and we conclude that $(s,t)$ is a subduction.
\end{proof}

\begin{remark}
The equivalence noted in \cref{p:gerbiness} is proven in the context of stacks over topological spaces in \cite[Proposition 5.3]{noohi}, albeit with respect to a more general definition of a gerbe.
\end{remark}

\section{Basic Differential Forms}\labell{s:differential forms}

Recall that given a Lie groupoid $G$ with source and target map $s$ and $t$, resp., a ``basic differential form'' $\alpha$ on $G_0$ is one such that $s^*\alpha=t^*\alpha$ (see \cite[Definition 3.1]{watts-groupoids}).  In this section we give a notion of a basic differential form on a stack, and show that it  extends the definition for a Lie groupoid (see \cref{p:differential forms}).  In the case of a stack equal to the Grothendieck construction applied to a concrete sheaf, the definition matches that for a corresponding diffeological space (see \cref{r:differential forms}).  It turns out that this definition only depends on the discretization of the stack, and in the case that the stack is differentiable and represented by a specific type of Lie groupoid, it only depends on the underlying diffeological coarse moduli space (see \cref{t:differential forms}).  It remains an open question whether this dependence on the diffeological coarse moduli space holds for all stacks.

\begin{definition}[Differential Forms]\labell{d:differential forms1}
Let $\Omega^k\colon\mfld^{op} \to \set$ be the sheaf of differential $k$-forms over manifolds.  Then, for a sheaf $S\colon\mfld^{op}\to \set$, we will define a \emph{basic differential $k$-form} to be a map of sheaves $\mu\colon S\to \Omega^k$.  Define $\Omega^k(S):=\Hom_\sheaves(S,\Omega^k)$.  Similarly, for a stack $\mathcal{X}$, define a \emph{basic $k$-form} to be a map of stacks $\beta\colon\mathcal{X}\to \int \Omega^k$ and let $\Omega^k(\mathcal{X}):=\Hom_{\stacks}(\mathcal{X},\int \Omega^k)$.
\end{definition}

\begin{remark}
Note that $\Omega^k(\mathcal{X})$ is discrete; \emph{i.e.}\ for each manifold $M$,  $\Omega^k(\mathcal{X})(M)$ is a set.  In fact $\Omega^k(\mathcal{X})(M)$ has a vector space structure, but we do not need this for our purposes.
\end{remark}

\begin{definition}[Differential Forms in Diffeology, see Article 6.28 of \cite{iglesias}]\labell{d:differential forms2}
A \emph{differential $k$-form} $\mu$ on a diffeological space $X$ is an assignment to each plot $p\colon U\to X$ a differential form $\mu(p)\in\Omega^k(U)$ satisfying the following compatibility condition: if $f\colon V\to U$ is a smooth map between open subsets of Euclidean spaces, then $\mu(p\circ f)=f^*\mu(p)$.
\end{definition}

\begin{remark}\labell{r:differential forms}
When considering a diffeological space $(X,\mathcal{D})$ as a concrete sheaf $I(\mathcal{D})\colon\mfld^{op} \to \set$, this definition says precisely that a differential $k$-form for a diffeological space is a map of sheaves $\mu\colon I(\mathcal{D}) \to \Omega^k$, matching \cref{d:differential forms1}.
\end{remark}

\begin{proposition}[Basic Forms of a Lie Groupoid]\labell{p:differential forms}
Let $G=(G_1\toto G_0)$ be a Lie groupoid, and let $y$ be the Yoneda embedding of $\mfld$ into $\stacks$.  An atlas $\alpha\colon y(G_0)\to BG$ induces a bijection from $\Omega^k(BG)$ to basic $k$-forms on $G_0$; that is, differential $k$-forms $\mu$ such that $s^*\mu=t^*\mu$ where $s$ and $t$ are the source and target maps of $G$, respectively.
\end{proposition}

\begin{proof}
Let $\mu\in\Omega^k(BG)$, and let $\alpha\colon y(G_0)\to BG$ be an atlas.  Then $$\mu\circ\alpha\circ y(s)=\mu\circ\alpha\circ y(t).$$  Thus, $\mu\circ\alpha$ corresponds to a basic differential form on $G_0$.

In the other direction, any basic form $\eta$ on $G_0$ induces a map of stacks $\tilde{\eta}\colon y(G_0)\to\int\Omega^k$ such that $$\tilde{\eta}\circ y(s)=\tilde{\eta}\circ y(t).$$  By \cref{l:BG 2-colimit} and the $2$-Yoneda Lemma, there is a map $\theta\colon BG\to\int\Omega^k$ such that $\tilde{\eta}=\theta\circ\alpha$.  Since $\int\Omega^k$ is discrete, $\theta$ is unique.
\end{proof}

\begin{theorem}[Differential Forms and $\coarse$]\labell{t:differential forms}
Let $\mathcal{X}$ be a stack.  Then, pre-composition with $\discr$ induces a bijection from $\Omega^k(\int\discr(\mathcal{X}))$ to $\Omega^k(\mathcal{X})$, and pre-composition with $\coarse$ induces an injection from $\Omega^k(\int\coarse(\mathcal{X}))$ to $\Omega^k(\mathcal{X})$.  Moreover, if $G$ is
\begin{enumerate}
\item \labell{i:groupoid} a finite-dimensional, Hausdorff, paracompact, proper Lie groupoid; or
\item \labell{i:gp action} the action groupoid of a Lie group on a manifold whose identity component acts properly; or
\item \labell{i:gerby} a gerbe over $G_0/G_1$,
\end{enumerate}
then the injection from $\Omega^k(\int\coarse(BG))$ to $\Omega^k(BG)$ is in fact a bijection.
\end{theorem}

\begin{proof}
The first claim is immediate from the fact that the Grothendieck construction is a fully faithful functor, along with \cref{l:discr adj}.  The second claim follows from the fact that any differential form $\mu\colon\coarse(\mathcal{X})\to\Omega^k$ is completely determined by it values on elements of the form $\underline{[\xi]}$ where $\xi$ are objects in $\mathcal{X}$.  Finally, the last claim follows from \cref{l:BH1}; \cref{r:differential forms}; and \cite{watts-groupoids} in the case of \cref{i:groupoid}, or \cite{KW} in the case of \cref{i:gp action}, or \cref{p:gerbiness} in the case of \cref{i:gerby}.
\end{proof}


\end{document}